\documentclass[english]{sig-alternate}
\usepackage[T1]{fontenc}
\usepackage[latin9]{inputenc}
\usepackage{microtype}
\usepackage{float}
\usepackage{wrapfig}
\usepackage{graphicx}
\usepackage[margin=0mm]{caption}
\usepackage{babel}
\usepackage{hyperref}
\usepackage{url}
\usepackage{algpseudocode}
\usepackage{algorithm}
\usepackage{color}

\definecolor{CBOrange}{RGB}{230,159,0}
\definecolor{CBBlue}{RGB}{86,180,233}
\definecolor{CBGreen}{RGB}{0,158,115}

\graphicspath{ {./img/} }

\def\sharedaffiliation{
\end{tabular}
\begin{tabular}{c}}

\newcommand*\rfrac[2]{{}^{#1}\!/_{#2}}

\makeatletter
\def\@copyrightspace{\relax}
\makeatother

\date{}
\hypersetup{colorlinks=true}
\begin{document}

\title{Wanted: Floating-Point Add Round-off Error instruction}

\numberofauthors{3}
\author{
\alignauthor
Marat Dukhan \\
	\email{mdukhan3@gatech.edu}
\alignauthor
Richard Vuduc
	\email{richie@cc.gatech.edu}
\alignauthor
Jason Riedy
	\email{jason.riedy@cc.gatech.edu}
\sharedaffiliation
	\affaddr{School of Computational Science and Engineering} \\
	\affaddr{College of Computing} \\
	\affaddr{Georgia Institute of Technology} \\
	\affaddr{Atlanta, GA}}

\maketitle
\begin{abstract}
We propose a new instruction (FPADDRE) that computes the round-off error in floating-point addition. We explain how this instruction benefits high-precision arithmetic operations in applications where double precision is not sufficient. Performance estimates on Intel Haswell, Intel Skylake, and AMD Steamroller processors, as well as Intel Knights Corner co-processor, demonstrate that such an instruction would improve the latency of double-double addition by up to $55\%$ and increase double-double addition throughput by up to $103\%$, with smaller, but non-negligible benefits for double-double multiplication. The new instruction delivers up to $2 \times$ speedups on three benchmarks that use high-precision floating-point arithmetic: double-double matrix-matrix multiplication, compensated dot product, and polynomial evaluation via the compensated Horner scheme.
\end{abstract}

\section{Introduction}

\begin{figure}[t]
	\centering
		\includegraphics[width=0.5\textwidth]{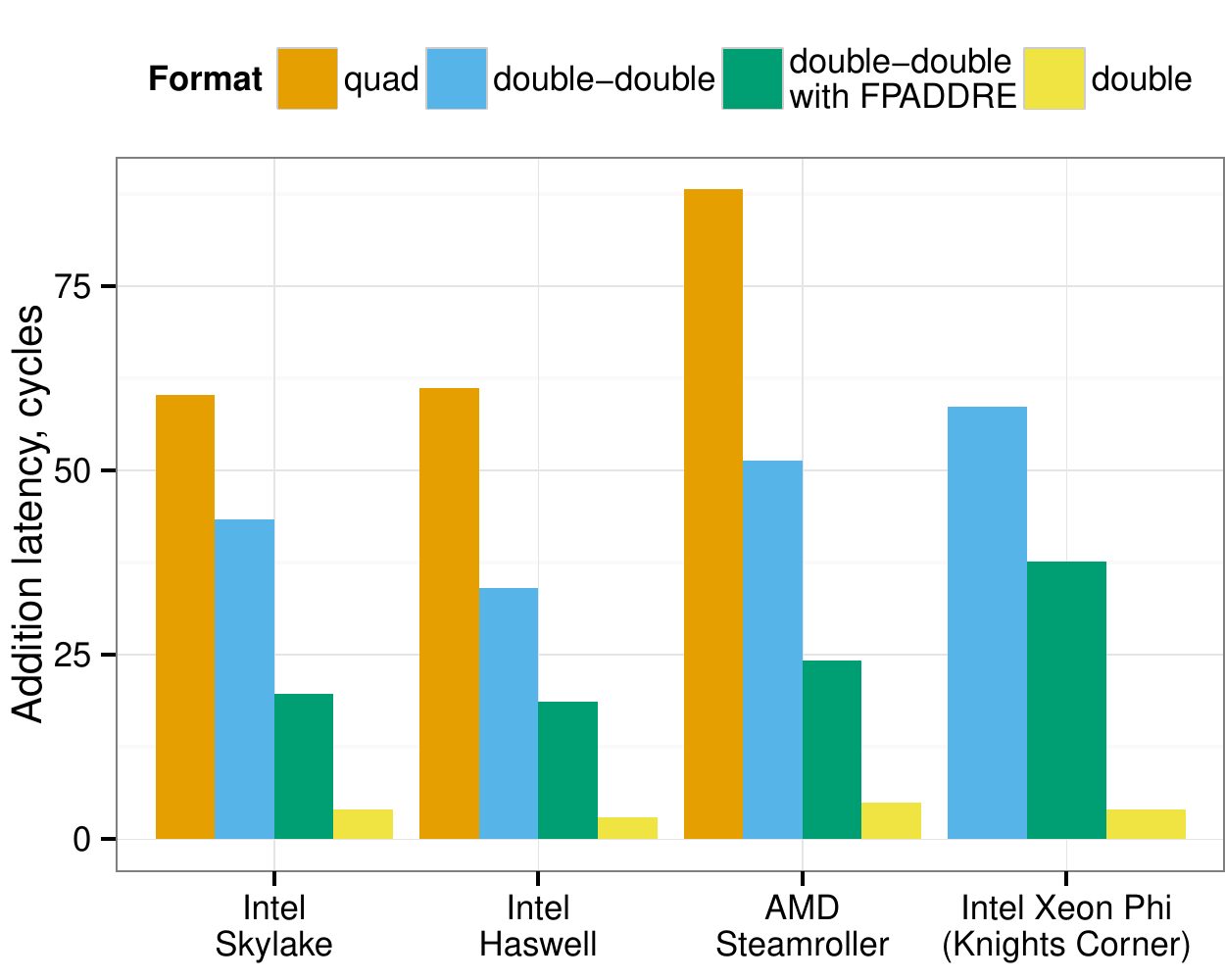}
	\caption{\label{fig:add-latency}Latency of floating-point addition operation with double-precision, double-double, and quad-precision floating-point formats}
\end{figure}

\textbf{High-precision floating-point computations} are represented by three kinds of algorithms: packed quadruple precision arithmetic, multi-word arithmetics like double-double, and compensated algorithms. Quadruple precision arithmetic is implemented either in rare hardware~\cite{Z990FP} or through integer operations in software (the \textcolor{CBOrange}{orange bars} in Fig.~\ref{fig:add-latency}). Double-double arithmetic (the \textcolor{CBBlue}{blue bars} in Fig.~\ref{fig:add-latency}) is implemented in software and extends precision, but not range, by representing a number as the unevaluated sum of a pair of double-precision values (or more for triple-double, etc.). Each double-double operation uses multiple double-precision operations to evaluate and renormalize the result. Compensated algorithms essentially inline the double-double operations and remove unnecessary intermediate normalizations for performance. This research advocates for a new instruction to greatly optimize (the \textcolor{CBGreen}{green bars} in Fig.~\ref{fig:add-latency}) the latter variants, double-double arithmetic, and compensated algorithms.

These arithmetics have been known for decades, but until recently they remained a rarely used hack of interest mostly to researchers in numerical computing. However, two recent trends suggest that high-precision floating-point arithmetics will become mainstream in the next decade. First, the wide availability of fused multiply-add (FMA) instructions on general-purpose hardware tremendously improves the performance of double-double operations and compensated algorithms. Secondly, new application requirements are making extended precision important for a wider audience:

\begin{itemize}
	\item \textbf{Numerical reproducibility} is an important issue on modern systems. First, with multiple separately clocked cores and non-uniform memory access it becomes inefficient to statically distribute work across many threads. Nondeterministic thread scheduling techniques, like work-stealing, help to exploit all available thread-level parallelism, but they also makes the result of floating-point summations dependent on random scheduling events, and thus non-reproducible. Secondly, variations in SIMD- and instruction-level parallelism across CPU and GPU architectures introduces similar numerical reproducibility issues across different platforms.  Computing intermediate results to extra precision can reduce reproducibility issues, and a slight modification to inner loops provides portable, accurate, and reproducible linear algebra based on double-double arithmetic~\cite{6875899}.
	\item The 2008 revision of IEEE-754 floating-point arithmetic standard recommends that \textbf{mathematical functions}, such as logarithms or arcsine, should be correctly rounded, i.e. accurate to the last bit. To satisfy this accuracy requirement, implementations need to use high-precision computations internally.
	\item The number of \textbf{scientific computing} applications that need more than double-precision arithmetic is increasing. David Bailey's review of high-precision floating-point arithmetic from 2005 lists 8 areas of science that use high-precision arithmetic~\cite{Bailey2005}, whereas his 2014 presentation has expanded the list to 12 areas~\cite{Bailey2014}.
\end{itemize}

\section{Error-Free Transformations}

\textbf{Error-free transformations} are the workhorses of both double-double arithmetic and compensated algorithms. Error-free addition represents the sum of two floating-point values $a + b$ as $s + e$ where $s$ is the result of floating-point addition instruction and $e$ is its round-off error. Similarly, error-free multiplication represents the product of two floating-point values $a \cdot b$ as $p + e$ where $p$ is the result of floating-point multiplication instructions and $e$ is the multiplication round-off error.

The multiplication round-off error can be computed with only one FMA instruction. However, computing the round-off error of addition in general requires 5 floating-point addition or subtraction instructions. In a special case when operands are ordered by magnitude, the round-off error can be computed with 2 floating-point instructions. Algorithms~\ref{alg:efadd-general} and~\ref{alg:efadd-special} illustrate the operations in the error-free addition for the general and the special cases.

\begin{algorithm}[h]
	\begin{algorithmic}
		\Function {Error-Free-Add-General}{$a$, $b$}
			\State $sum \gets \Call{FPADD}{a, b}$
			\State $b_{virtual} \gets \Call{FPADD}{sum, -a}$
			\State $a_{virtual} \gets \Call{FPADD}{sum, -b_{virtual}}$
			\State $b_{roundoff} \gets \Call{FPADD}{b, -b_{virtual}}$
			\State $a_{roundoff} \gets \Call{FPADD}{a, -a_{virtual}}$
			\State $error \gets \Call{FPADD}{a_{roundoff}, b_{roundoff}}$
			\State \textbf{return} $sum, error$
		\EndFunction
	\end{algorithmic}
	\caption{\label{alg:efadd-general}Error-free addition algorithm for the general case. The algorithm is due to Knuth~\cite{Knuth1997}, but the listing below follows the notation of Theorem 7 from Shewchuk~\cite{Shewchuk1997}.}
\end{algorithm}

\begin{algorithm}[h]
	\begin{algorithmic}
		\Function {Error-Free-Add-Special}{$a$, $b$}
			\State $sum \gets \Call{FPADD}{a, b}$
			\State $b_{virtual} \gets \Call{FPADD}{sum, -a}$
			\State $error \gets \Call{FPADD}{b, -b_{virtual}}$
			\State \textbf{return} $sum, error$
		\EndFunction
	\end{algorithmic}
	\caption{\label{alg:efadd-special}Error-free addition algorithm for the special case when $|a| \geq |b|$. The algorithm is due to Dekker~\cite{Dekker1971}, but the listing below follows the notation of Theorem 6 from Shewchuk~\cite{Shewchuk1997}.}
\end{algorithm}

\subsection{FPADDRE Instruction}

\begin{algorithm}
	\begin{algorithmic}
		\Function {Error-Free-Add-With-FPADDRE}{$a$, $b$}
			\State $sum \gets \Call{FPADD}{a, b}$
			\State $error \gets \Call{FPADDRE}{a, b}$
			\State \textbf{return} $sum, error$
		\EndFunction
	\end{algorithmic}
	\caption{\label{alg:efadd-fpaddre}Error-free addition algorithm in the general case using the proposed \textbf{FPADDRE} instruction. Note that the two operations in the algorithm are independent of each other, and could be computed in parallel.}
\end{algorithm}

We propose a new instruction, \textbf{Floating-Point Addition Round-off Error (FPADDRE)}, that complements floating-point addition instruction (FPADD), and makes possible to compute error-free addition in just two instructions, as demonstrated in Alg.~\ref{alg:efadd-fpaddre}.

\begin{figure}[h]
	\includegraphics[width=\hsize]{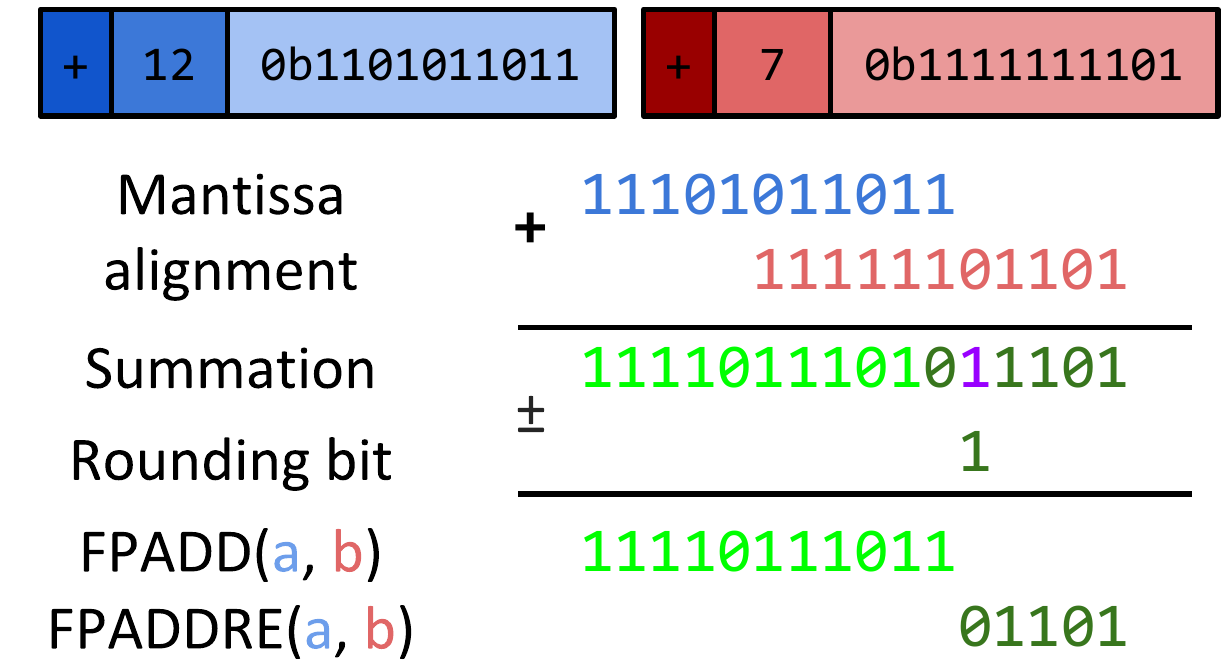}
	\caption{\label{fig:fpaddre-schema}Schema of FPADD and FPADDRE operations (the case of operands with the same sign and overlapping mantissas). The operations differ only in two aspects: addition or subtraction of a sticky bit and the bits copied to the resulting mantissa.}
\end{figure}

The floating-point addition (FPADD) instruction computes the sum of two floating-point numbers and then rounds the result to the nearest floating-point number, losing information in the last bits of the sum. The proposed FPADDRE instruction performs a similar operation but returns the last, normally wasted, bits of the sum. Figure~\ref{fig:fpaddre-schema} illustrates the similarities and differencies of the two operations. FPADDRE differs only slightly from addition and could reuse its circuits in a hardware implementation. Besides replacing $5$ FPADD operations, FPADDRE additionally improves the latency of error-free transformation by breaking the dependency chain between the addition result and the round-off error.

\begin{algorithm}
	\begin{algorithmic}
		\Function {Error-Free-Add-With-FPADD3}{$a$, $b$}
			\State $sum \gets \Call{FPADD}{a, b}$
			\State $error \gets \Call{FPADD3}{a, b, -sum}$
			\State \textbf{return} $sum, error$
		\EndFunction
	\end{algorithmic}
	\caption{\label{alg:efadd-fpadd3}Error-free addition algorithm in the general case using the \textbf{FPADD3} instruction, suggested by Ogita et al.~\cite{Ogita2005} Note that the two operations in the algorithm form a dependency chain and \textit{cannot} be computed in parallel.}
\end{algorithm}

Ogita et al. proposed the FPADD3 instruction, which adds three floating-point values without intermediate rounding~\cite{Ogita2005}. Algorithm~\ref{alg:efadd-fpadd3} shows that with an FPADD3 instruction, it would be possible to compute the addition round-off error with one instruction, albeit with a dependency on the FPADD result. Unfortunately, no general-purpose CPU or GPU implements the FPADD3 operation. One reason may be that a fast hardware implementation of FPADD3 would require considering 4 overlapping options for the three inputs, compared to just 2 overlapping cases in regular addition. The suggested FPADDRE instruction does not share this drawback.

\begin{table}[!h]
	\begin{tabular}{ l | c | c }
		Error-Free Addition           & Instructions & Latency slots \\
		\hline
		General case                  & 6            & 5             \\
		Special case ($|a| \geq |b|$) & 3            & 3             \\
		With FPADDRE                  & 2            & 1             \\
		With FPADD3                   & 2            & 2             \\
	\end{tabular}
	\caption{\label{tab:efadd-perf}Performance characteristics of error-free addition algorithm in different implementations.}
\end{table}

Table~\ref{tab:efadd-perf} summarizes the performance characteristics of the four versions of an error-free addition algorithm. It shows that an FPADDRE instruction enables the most performant implementation.

\section{Performance simulation}

We evaluate the speedups achievable with hardware FPADDRE implementations on three recent x86-64 processor microarchitectures from Intel and AMD as well as on the Intel Xeon Phi co-processor based on the Knights Corner microarchitecture. Table~\ref{tab:processors} details the benchmarking platforms. We evaluated all processors in single-thread mode, which can be suboptimal for absolute performance; however, we have no reasons to expect that it leads to systematic errors in estimation of speedups due to FPADDRE instruction. Because FPADDRE is similar to floating-point addition, we assume that a hardware implementation would exhibit the same performance characteristics as floating-point addition. We implemented several high-precision floating-point benchmarks in C with intrinsics and ran two sets of tests. In the first set, all operations were implemented with the default instruction set of the respective architectures. In the second set of tests, we simulated a FPADDRE instruction by replacing it with an instruction that has the same performance characteristics as floating-point addition. On AMD Steamroller we simulated \texttt{fpaddre(a, b)} as \texttt{fma(a, a, b)} and on other architectures we replaced it with \texttt{min(a, b)}. Of course, such substitutions may lead to incorrect numerical results, but the incorrect results do not affect the control flow of the benchmarks.

\begin{table*}[f]
	\begin{tabular}{ l | c | c | c | c }
		Processor               & Intel Core i7-6700K & Intel Core i7-4700K & AMD A10-7850K & Intel Xeon Phi SE10P          \\
		\hline
		Microarchitecture       & Skylake             & Haswell             & Steamroller       & Knights Corner            \\
		Frequency               & $4.0$ GHz           & $3.5$ GHz           & $3.7$ GHz         & $1.1$ GHz                 \\
		\hline
		L1D Cache               & $32K$               & $32K$               & $16K$             & $32K$                     \\
		L2 Cache                & $256K$              & $256K$              & $2M$              & $512K$                    \\
		L3 Cache                & $8M$                & $8M$                & None              & None                      \\
		\hline
		SIMD width (double)     & $4$                 & $4$                 & $4$               & $8$                       \\
		SIMD ADD issue ports    & P0 or P1            & P0                  & P0 and P1         & VALU                      \\
		SIMD MUL issue ports    & P0 or P1            & P0 or P1            & P0 and P1         & VALU                      \\
		SIMD FMA issue ports    & P0 or P1            & P0 or P1            & P0 and P1         & VALU                      \\
		FP ADD latency          & $4$                 & $3$                 & $5$               & $4$                       \\
		FP MUL latency          & $4$                 & $5$                 & $5$               & $4$                       \\
		FP FMA latency          & $4$                 & $5$                 & $5$               & $4$                       \\
		\hline
		Compiler                & gcc 5.2.1           & gcc 5.2.1           & gcc 5.2.1         & icc 15.0.0                \\
		Optimization flags      & -O3 -mavx2 -mfma    & -O3 -mcore-avx2     & -O3 -march=bdver3 & -O3 -mmic                 \\
		Floating-point flags    & -ffp-contract=off   & -ffp-contract=off   & -ffp-contract=off & -fp-model precise -no-fma \\
	\end{tabular}
	\caption{\label{tab:processors}Processors and co-processor used in performance evaluation}
\end{table*}

\subsection{Microbenchmarks}

In this set of benchmarks, we measured the effect of FPADDRE instruction on the performance characteristics of double-double addition and multiplication operations. The double-double addition involves 2 general-case error-free additions and 2 special-case error-free additions, and benefits from FPADDRE the most. The double-double multiplication involves only one special-case error-free addition, but nonetheless exhibits some speedup from FPADDRE.

Each microbenchmark was replicated $1000$ times, with the best performance reported. We observed negligible variation in reported performance across independent runs.

\begin{figure}[h]
	\centering
		\includegraphics[width=\hsize]{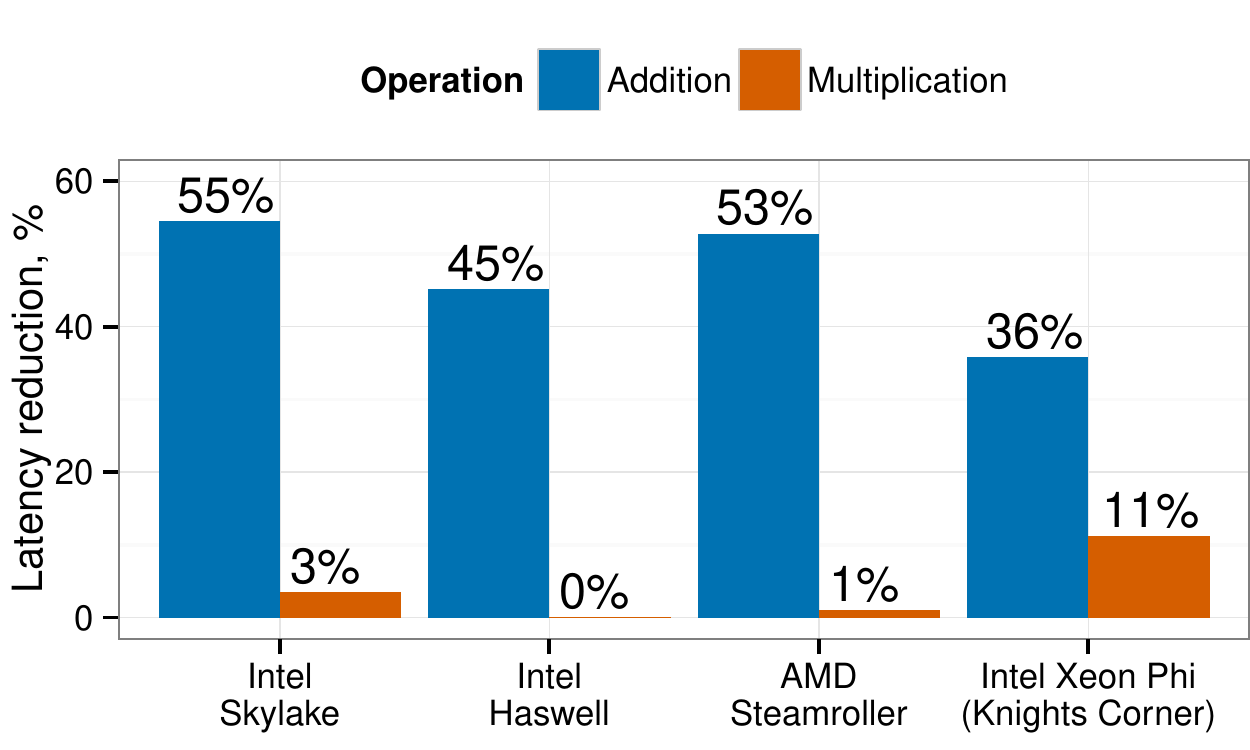}
	\caption{\label{fig:doubledouble-latency}Reduction of latency of double-double addition and multiplication due to FPADDRE instruction}
\end{figure}

Figure~\ref{fig:doubledouble-latency} shows the reduction of latency due to FPADDRE for double-double addition and multiplication. In this benchmark, we measured the time to add or multiply all elements of a double-double array that fits in the L1 cache.

\begin{figure}[h]
	\centering
		\includegraphics[width=\hsize]{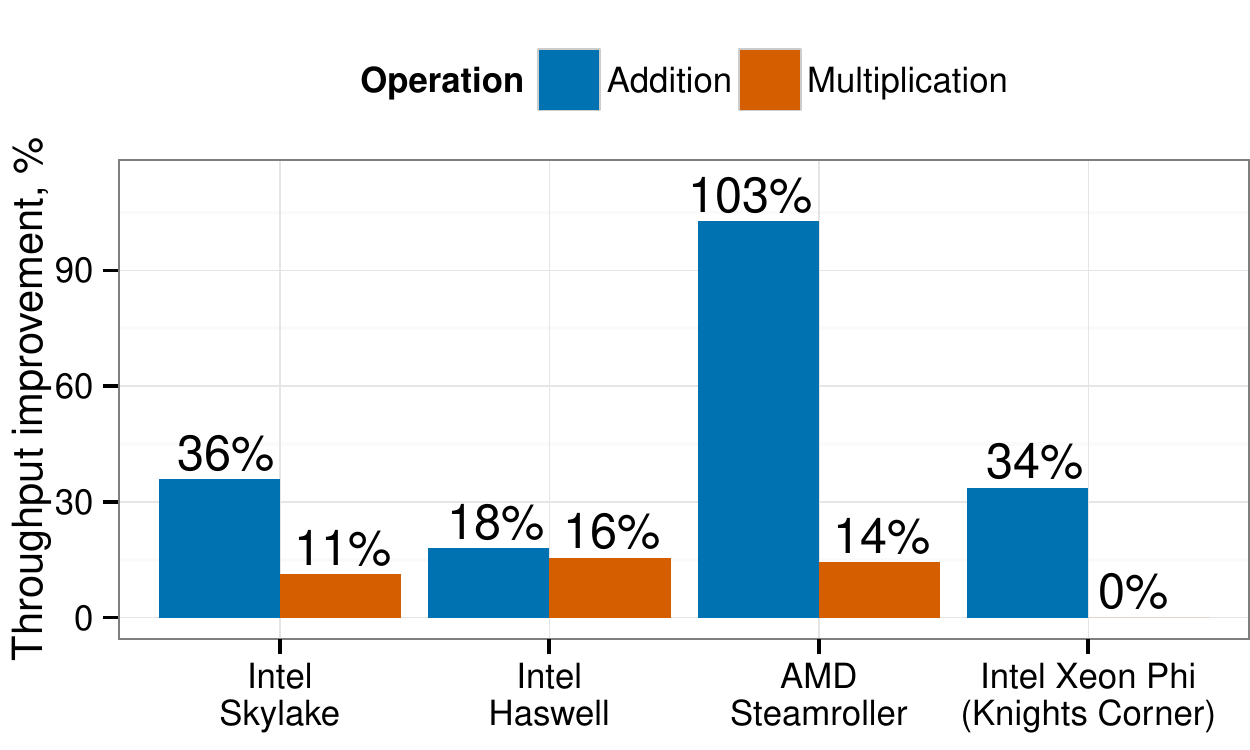}
	\caption{\label{fig:doubledouble-throughput}Throughput improvement of double-double addition and multiplication due to FPADDRE instruction}
\end{figure}

Figure~\ref{fig:doubledouble-throughput} shows how FPADDRE instruction improves throughput of double-double addition and multiplication. In this benchmark, we measured the time to add or multipliply a double-double array with a double-double constant. In this case, the operations on different array elements are independent and can be performed in parallel. The array size was selected to fit into the L1 cache.

\subsection{Applications}

Beyond low-level microbenchmarks, we also profiled three kernels that arise in important applications of high-precision arithmetic: polynomial evaluation with the compensated Horner scheme, compensated dot product, and double-double matrix multiplication. Each application benchmark was repeated at least $1000$ times, and we report the median performance across runs.

\begin{figure}[h]
	\centering
		\includegraphics[width=\hsize]{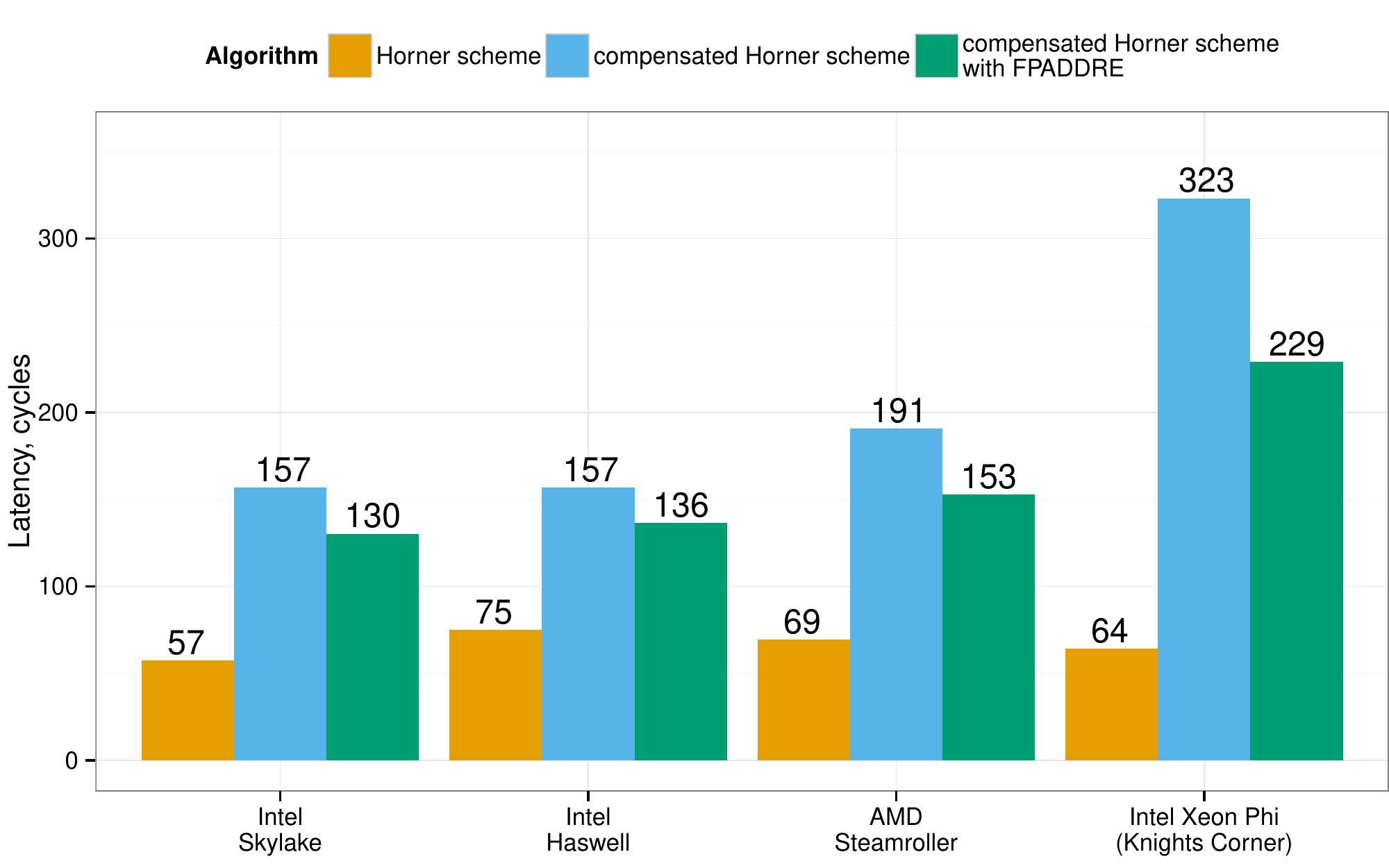}
	\caption{\label{fig:polevl-latency}Latency of 15-order polynomial evaluation with compensated Horner scheme}
\end{figure}

The compensated Horner scheme evaluates polynomial with double-precision coefficients in approximately double-double intermediate precision~\cite{CompensatedHorner2005}. This algorithm is useful in correctly-rounded implementations of mathematical functions. Figure~\ref{fig:polevl-latency} shows that a FPADDRE instruction reduces the latency of 15-degree polynomial evaluation by $13\% - 29\%$ for the microarchitectures considered.

\begin{figure*}[t]
	\centering
		\includegraphics[width=\hsize]{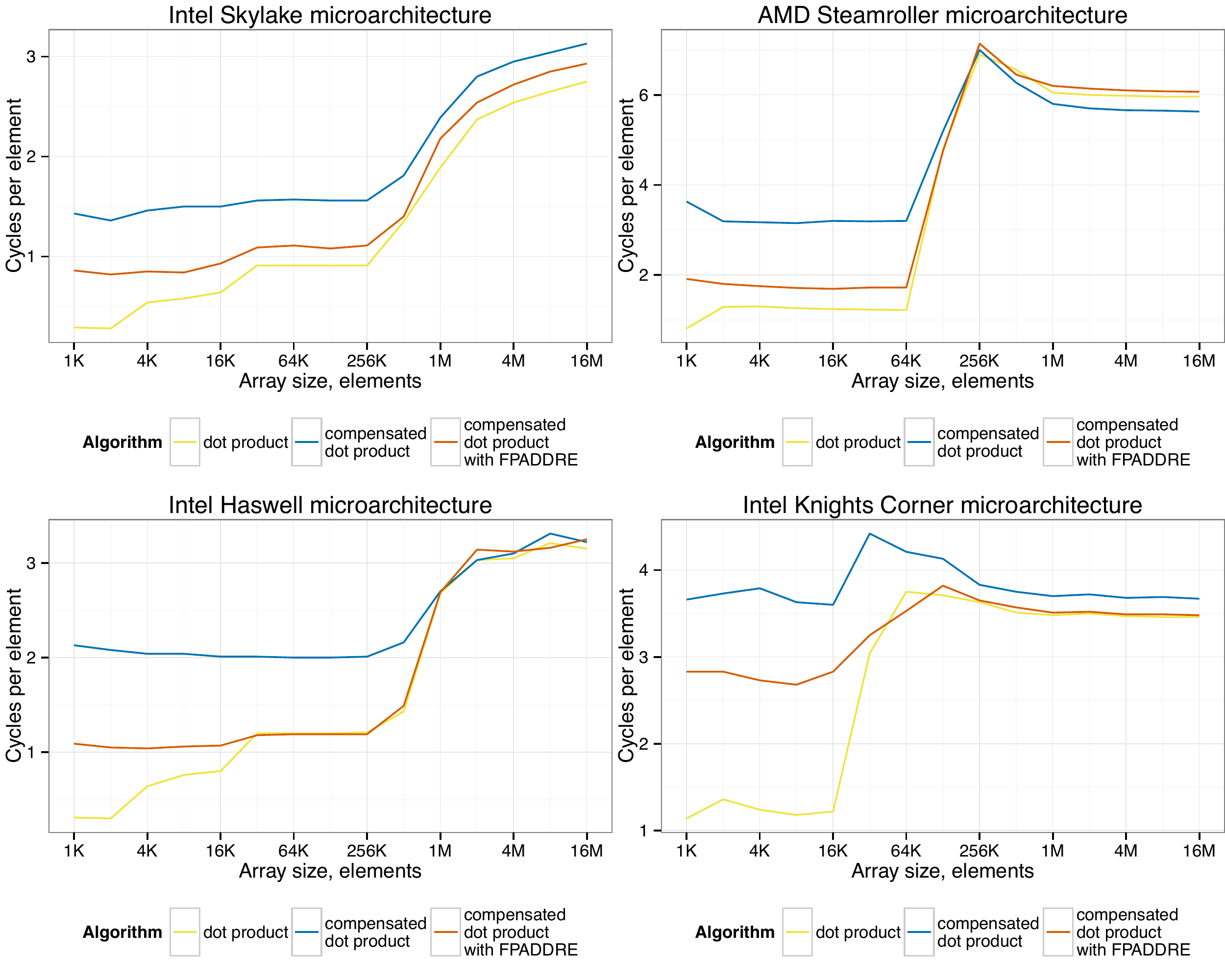}
	\caption{\label{fig:dot-product}Dot product and compensated dot product performance with standard ISA and with FPADDRE instruction}
\end{figure*}

The compensated dot product algorithm computes the dot product of double-precision vectors in approximately double-double internal precision~\cite{Ogita2005}. The extra precision helps reproducibility on large data sets. For the benchmark, we implemented dot product and compensated dot product algorithms using SIMD intrinsics and unrolled the main loop by factors of $1$ to $8$. Figure~\ref{fig:dot-product} shows the performance of the dot product algorithms with the most performant unroll factors for each microarchitecture, algorithm, and array size. When the arrays fit into the L1 cache, the compensated dot product algorithm is compute-bound, and the FPADDRE instruction increases performance by $66\%$, $95\%$, $93\%$, and $29\%$ on Intel Skylake, Intel Haswell, AMD Steamroller, and Intel Knights Corner, respectively.
One additional instruction per operand guarantees reproducibility in summation~\cite{6875899}, so fully reproducible dot products should see similar performance improvements.

\begin{figure}[h]
	\centering
		\includegraphics[width=\hsize]{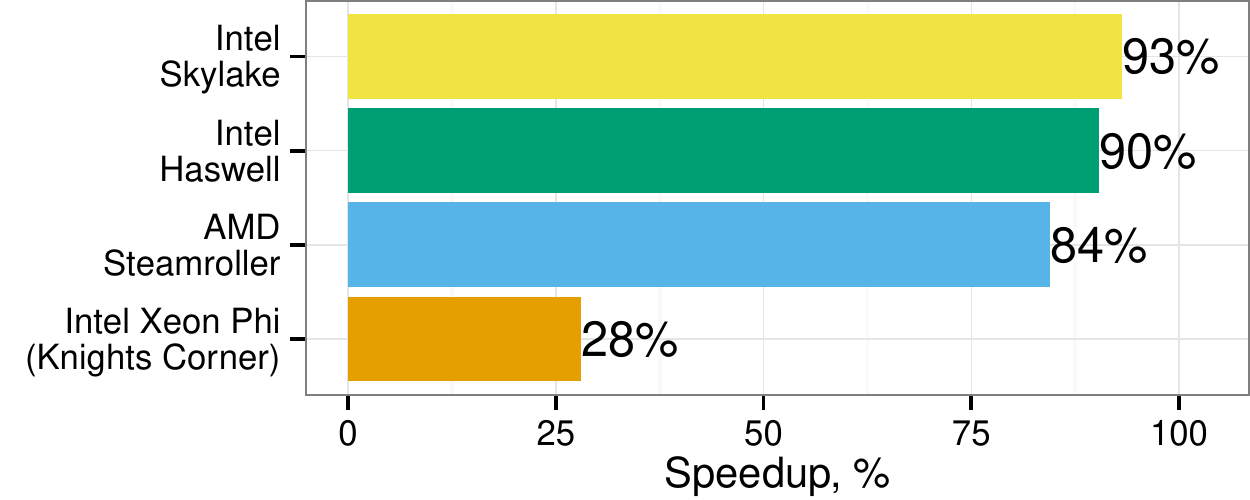}
	\caption{\label{fig:gemm-speedup}Double-double matrix multiplication acceleration with FPADDRE instruction}
\end{figure}

Matrix-matrix multiplication is a basic building block for computational linear algebra algorithms. Recently, van Zee and van de Geijn showed that state-of-the-art high-performance implementations of matrix-matrix multiplication can be written mostly in portable high-level code, with only the small inner kernels using target-specific intrinsics or assembly~\cite{BLIS}. For another benchmark, we implemented the inner kernel of general matrix-matrix multiplication (DDGEMM), where inputs, outputs and intermediate values are stored in double-double format. Typically, this inner kernel is responsible for over $90\%$ of compute time in a matrix-matrix multiplication. We considered inner kernels with multiple register blocking parameters and selected the parameters that deliver the best performance for each microarchitecture. Figure~\ref{fig:gemm-speedup} characterises the speedups in the matrix multiplication when the instruction set is enriched with FPADDRE, and Table~\ref{tab:ddgemm-perf} characterises absolute performance of double-double matrix multipication micro-kernel (DDGEMM) in double-double MFLOPS and compares it to double-precision matrix multiplication (DGEMM) in double-precision MFLOPS. DGEMM performance is measured on production-quality libraries with $4096\times4096$ matrices; the DDGEMM performance numbers are for the micro-kernel only and ignore the overhead of repacking the matrices and boundary effects. The data in Table~\ref{tab:ddgemm-perf} demonstrates that double-double matrix multiplication is presently $35-69 \times$ slower than in double-precision, and thus FPADDRE-provided acceleration is very welcome.

\begin{table*}[f]
	\begin{tabular}{ l | c | c | c | c }
		Operation           & Intel Skylake                    & Intel Haswell                    & AMD Steamroller                  & Intel Knights Corner            \\
		\hline
		DDGEMM              & $1732~(\approx\rfrac{1}{37}$ DP) & $1199~(\approx\rfrac{1}{45}$ DP) &  $743~(\approx\rfrac{1}{25}$ DP) & $255~(\approx\rfrac{1}{17}$ DP) \\
		DDGEMM with FPADDRE & $3344~(\approx\rfrac{1}{19}$ DP) & $2283~(\approx\rfrac{1}{24}$ DP) & $1370~(\approx\rfrac{1}{19}$ DP) & $326~(\approx\rfrac{1}{14}$ DP) \\
		\hline
		DGEMM               & $63603$ (MKL)                    & $51409$ (MKL)                    & $25869$ (OpenBLAS)               & $4439$ (MKL)                    \\
	\end{tabular}
	\caption{\label{tab:ddgemm-perf}Performance (in MFLOPS) of general matrix-matrix multiplication on the four benchmarked microarchitectures}
\end{table*}

\section{FPMULRE Instruction}

\begin{algorithm}
	\begin{algorithmic}
		\Function {Error-Free-Mul-With-FMA}{$a$, $b$}
			\State $product \gets \Call{FPMUL}{a, b}$
			\State $error \gets \Call{FMA}{a, b, -product}$
			\State \textbf{return} $product, error$
		\EndFunction
	\end{algorithmic}
	\caption{\label{alg:efmul-fma}Error-free multiplication algorithm using the \textbf{FMA} instruction. The two operations in the algorithm form a dependency chain and \textit{cannot} be computed in parallel.}
\end{algorithm}

An error-free multiplication transformation represents the product of two floating-point values $a \cdot b$ as $p + e$, where $p$ is the result of the floating-point multiplication instruction and $e$ is the multiplication round-off error. Most modern hardware platforms implement fused multiply-add (FMA) instructions, which permit computation of error-free multiplication with just two instructions, as illustrated in Alg.~\ref{alg:efmul-fma}.

\begin{algorithm}
	\begin{algorithmic}
		\Function {Error-Free-Mul-With-FPMULRE}{$a$, $b$}
			\State $product \gets \Call{FPMUL}{a, b}$
			\State $error \gets \Call{FPMULRE}{a, b}$
			\State \textbf{return} $product, error$
		\EndFunction
	\end{algorithmic}
	\caption{\label{alg:efmul-fpmulre}Error-free multiplication algorithm using the proposed \textbf{FPMULRE} instruction. Note that the two operations in the algorithm are independent of each other and can be computed in parallel.}
\end{algorithm}

Nonetheless, a case can be made for an instruction similar to FPADDRE, but for multiplication. We would call this instruction Floating-Point Multiplication Round-off Error (FPMULRE). FPMULRE would directly compute the round-off error of multiplication operations; the operation is similar to floating-point multiplication and could share hardware circuits with it. FPMULRE would benefit error-free multiplications in two ways. First, it would allow computation of both outputs of error-free multiplication in parallel, thereby reducing its latency. Algorithm~\ref{alg:efmul-fpmulre} illustrates this point. Secondly, FPMULRE is simpler than fused multiply-add, and its implementation can be more energy-efficient than an FMA's implementation. As FPMULRE replaces FMA in the error-free multiplication, it could result in better energy efficiency of this transformation and higher-level high-precision operations.

\section{Conclusions}

High-precision floating-point arithmetic is about to become a common computing technique. This claim motivates our proposed Floating-Point Addition Round-off Error (FPADDRE), which can accelerate error-free addition and algorithms based on error-free transformations, including double-double arithmetics, compensated Horner scheme, compensated dot product, and compensated summation. Our performance simulations suggest that if FPADDRE were implemented in recent processors and co-processors, we would observe a $13\%-29\%$ reduction in latency of compensated Horner scheme, up to $29\%-95\%$ performance increase in compensated dot product, and $28\%-93\%$ speedup in double-double matrix multiplication. The same idea could be translated to multiplication, where a Floating-Point Multiplication Round-off Error (FPMULRE) would improve latency and energy efficiency of error-free multiplication.

To facilitate and encourage further research on these ideas we released the source code for the implemented algorithms and simulations on \href{https://github.com/Maratyszcza/FPplus}{www.GitHub.com/Maratyszcza/FPplus}.

\section*{Acknowledgements}

We thank Edmond Chow and his Intel Parallel Computing Center at Georgia Tech for access to Xeon Phi-based platform.  We also thank James Demmel (UC Berkeley) and Greg Henry (Intel) for discussions on reproducibility and performance.

This material is based upon work supported by the U.S. National Science Foundation (NSF) Award Number 1339745 
and the U.S. Dept. of Energy (DOE), Office of Science, Advanced Scientific Computing Research under award DE-FC02-10ER26006/DE-SC0004915. 
Any opinions, findings and conclusions or recommendations expressed in this material are those of the authors and do not necessarily reflect those of NSF or DOE.
The work also partially sponsored by Defense Advanced Research Projects Agency (DARPA) under agreement \#HR0011-13-2-0001. The content, views and conclusions presented in this document do not necessarily reflect the position or the policy of DARPA or the U.S. Government, no official endorsement should be inferred. Distribution Statement A: ``Approved for public release; distribution is unlimited.'' 

\bibliography{fpadde}{}
\bibliographystyle{plain}

\end{document}